\documentclass{ifacconf}
\usepackage{amsmath, amssymb}

\counterwithin*{section}{part}
\usepackage{cleveref}

\usepackage{graphicx,subfigure}      
\usepackage{algorithm}
\usepackage{algorithmicx}
\usepackage{algpseudocode}
\usepackage[numbers]{natbib}        
\usepackage{tikz}
\usetikzlibrary{arrows.meta}
\usetikzlibrary{cd,calc}
\tikzset{>=latex}
\graphicspath{{./Figs/}}

\usepackage{dsfont}

\definecolor{myblue}{RGB}{87, 163, 213}
\definecolor{myorange}{RGB}{226, 112, 46}
\definecolor{myred}{RGB}{218, 0, 0}
\definecolor{mygreen}{RGB}{129, 214, 83}
\definecolor{mygray}{gray}{0.9}

\newcommand{\Input}[1]{\textbf{Input:} #1}

\algblockdefx{FnBlock}{EndFnBlock}%
	[1]{\textbf{def} #1}
	{\textbf{return}}

\newcommand{\norm}[1]{\ensuremath{\left\| #1 \right\|}}

\newcommand{\bracket}[1]{\ensuremath{\left[ #1 \right]}}
\newcommand{\braces}[1]{\ensuremath{\left\{ #1 \right\}}}
\newcommand{\parenth}[1]{\ensuremath{\left( #1 \right)}}
\newcommand{\E}{\mathbb{E}}

\newcommand{\pair}[1]{\ensuremath{\left\langle #1 \right\rangle}}
\renewcommand{\Re}{\ensuremath{\mathbb{R}}}

\newcommand{\deriv}[2]{\ensuremath{\frac{\partial #1}{\partial #2}}}

\newcommand{\diver}{\text{div}}

\newtheorem{definition}{Definition}

\newtheorem{remark}{Remark}

\begin{document}
\begin{frontmatter}

\title{Mean Field Game and Control for Switching Hybrid Systems\thanksref{footnoteinfo}} 
\thanks[footnoteinfo]{This research has been supported in part by AFOSR under the grant FA9550-23-1-0400.}

\author[First]{Tejaswi K. C.} 
\author[Second]{Taeyoung Lee}

\address[First]{Mechanical and Aerospace Engineering, George Washington University, DC 20052 (email: kctejaswi999@gwu.edu)}
\address[Second]{Mechanical and Aerospace Engineering, George Washington University, DC 20052 (email: tylee@gwu.edu)}

\begin{abstract}                
    Mean field games and controls involve guiding the behavior of large populations of interacting agents, where each individual’s influence on the group is negligible but collectively impacts overall dynamics.
    Hybrid systems integrate continuous dynamics with discrete transitions, effectively modeling the complex interplay between continuous flows and instantaneous jumps in a unified framework.
    This paper formulates mean field game and control strategies for switching hybrid systems and proposes computational methods to solve the resulting integro-partial differential equation.
    This approach enables scalable, decentralized decision-making in large-scale switching systems, which is illustrated through numerical examples in an emergency evacuation scenario with sudden changes in the surrounding environment.
\end{abstract}


\end{frontmatter}

\section{Introduction}

Mean field game (MFG) theory and mean field control (MFC) provide a powerful framework for modeling large-scale systems where individual agents interact in a decentralized manner but influence each other’s behavior through aggregated effects \citep{lasry2007mean,caines2021mean}.
They also offer methods for efficiently controlling and optimizing complex systems with massive populations of agents.
Traditional centralized control approaches may become computationally prohibitive or infeasible when dealing with large numbers of interconnected entities.
By considering the aggregate behavior of the population, MFG/MFC enables the design of decentralized control strategies that are scalable and computationally efficient.
This is crucial in engineering applications where real-time decision-making and coordination are essential, such as in smart grids, swarm robotics, and crowd management systems \citep{cousin2011mean,huang2003individual}.

Hybrid dynamical systems are characterized by the interaction of continuous and discrete dynamics, where the continuous component represents processes evolving over time according to differential equations, and the discrete component captures sudden transitions or events, such as switching between operational modes~\citep{goedel2012hybrid,van2007introduction}.
This dual nature allows hybrid systems to model scenarios where processes are driven both by smooth, continuous changes and by abrupt shifts, often resulting from external triggers or internal control decisions.
By unifying continuous and discrete dynamics, hybrid systems enable precise modeling and control of critical applications that must reliably respond to both gradual variations and sudden events, enhancing system robustness, efficiency, and safety.

Mean field games can be integrated with hybrid dynamical systems, where discrete events or transitions are triggered by external conditions, thresholds, or state-dependent events within interactions among large populations.
For example, this approach could model traffic flow with vehicles that switch between different driving modes (e.g., autonomous vs. manual control) based on traffic conditions and efficiency.
Hybrid mean field games have been formulated for linear Gaussian systems to investigate optimal switching and stopping times~\citep{firoozi2022class}, where each agent has stochastic linear dynamics with quadratic costs.
This framework is applied to optimal execution problems in financial markets~\citep{firoozi2017mean}.

This paper addresses mean field games and controls for a more general class of hybrid systems that are not limited to quadratic costs for linear dynamics.
Specifically, we consider cases where cost functions may include nonlinear, coupled terms between the mean field and the control input, and where the surrounding environment may undergo abrupt transitions in boundary conditions.
We focus on computational methods to solve the governing equations for the proposed mean field games and controls for switching hybrid systems, represented by coupled integro-partial differential equations involving switching.
Numerical solutions for these coupled forward-backward equations are obtained by discretizing space and time.

For a specific example of emergency evacuation, we investigate the unique behaviors in these hybrid settings and examine how they differ from those in standard mean field games.
Additionally, we compare this solution with the case of the social optimum via mean field controls, in which agents coordinate with one another rather than competing.

In summary, the main contributions of this paper are the formulation of mean field games and controls for a specific class of hybrid systems and a computational framework for solving the resulting coupled governing equations under sudden environmental changes.

This paper is organized as follows.
Mean field game theory with continuous dynamics is reviewed in \Cref{sec:MFG}, after which the numerical algorithm to solve the hybrid mean field game problem is presented in \Cref{sec:ALG}, followed by numerical examples in \Cref{sec:SIM}.

\section{Mean Field Game and Control}\label{sec:MFG}

In this section, we formulate the classical $ N $-agent games where the interchangeability of individual agents does not alter the game's properties.  
The concept of Nash equilibria, as discussed in \cite{cardaliaguet2018short}, is employed to characterize the optimal strategies. 
We then present the governing equations for the density and value function in the limit as $N \to \infty$, from which the control input can be determined. 
Specifically, the desired response of an individual agent is obtained by minimizing the corresponding Hamiltonian of the control problem, using the equilibrium condition of the distribution. 
However, the evolution of the population density is inherently coupled with these equations, as it is generated based on the optimal policies of the agents.  

\subsection{Finite agent games}

Let us first consider a system of $ N $-agents, where $ N $ is a potentially large integer.
Each agent's dynamics can be expressed as:
\begin{equation}\label{eqn:dx_i}
dx^i(t) = \frac{1}{N}\sum_{j=1}^{N} b^i(x^i(t), x^j(t), \alpha^i(t)) dt + \sigma d W^i(t),
\end{equation}
for each $ i = 1,\dots, N $.
Here, $ b^i : \Re^n \times \Re^n \times \Re^d \to \Re^n $ is the drift term, and $ \alpha^i(t) \in \Re^d $ is the control for each agent. 
The initial states $ x^i(0) $ are distributed i.i.d. in $ \Re^n $, and $ W^i $ are independent $ n $-dimensional Brownian motions.
The parameter $ \sigma \in \Re $ is the scale factor for the noise terms.

There are a few key assumptions in this scenario:
\begin{itemize}
	\item The system is \textit{homogeneous}, meaning that the dynamical properties and the cost functions have the same structure for all agents:
	\begin{gather*}
	b^i = b^j \quad \forall i, j
	\end{gather*}
	\item Also, there is \textit{anonymity}, implying that interactions happen only through the aggregated mean field :
	\begin{equation}\label{eqn:mN}
	m^N(t) = \frac{1}{N} \sum_{j=1}^{N} \delta_{x^j(t)},
	\end{equation}
	which is an estimate of the population distribution.
\end{itemize}

Initialized from $ x^i(0) \sim m^N(0) $, the dynamics in \eqref{eqn:dx_i} is simplified into
\begin{equation}\label{eqn:dx_im}
dx^i(t) = b(x^i(t), m^N(t), \alpha^i(t, x^i(t))) dt + \sigma d W^i(t),
\end{equation}
where $ \alpha^i(t, x^i(t)) $ is the closed-loop control.

Each agent aims to minimize the following cost by solving an optimization problem over $ [0, T] $ :
\begin{multline}\label{eqn:Ji}
J^i(\alpha^i) = \E\bracket{\int_{0}^{T} f(x^i(t), m^N(t), \alpha^i(t)) dt + \right. \\
	\left. g(x^i(T), m^N(T)) },
\end{multline}
where $ f : \Re^n \times \Re^n \times \Re^d \to \Re $ denotes the running cost, and $ g : \Re^n \times \Re^n \to \Re $ is the terminal cost, both of which are homogeneous across all agents. 

Let $ \hat{\alpha}^{-i} $ represent the strategies of the remaining agents. 
The optimal strategy $ \hat{\alpha} = (\hat{\alpha}^i, \hat{\alpha}^{-i}) $ for all agents is determined as the Nash equilibrium of this $ N $-agent game.
That is, the current agent's $ \hat{\alpha}^i $ minimizes its own cost
\begin{equation}\label{eqn:nashN}
J^i(\hat{\alpha}^i) \le J^i(\alpha^i), \quad \forall i, \forall \alpha^i.
\end{equation}
which is defined in \eqref{eqn:Ji}.

\subsection{Mean field game : large $ N $ limit}

In the limit of $ N \to \infty $, the aggregate distribution $ m^N(t) $ is identified with a population state $ m(t) \in L^2(\Re^n) $. 
Mean field games are employed to study the Nash equilibria of such systems, providing an approximation for the behavior of a large number of agents. 
It is sufficient to understand the behavior of one representative agent since each individual has a negligible impact on the rest of the population. 

From \eqref{eqn:dx_im}, the dynamics of a representative agent is reduced to
\begin{equation}\label{eqn:dx}
	dx(t) = b(x(t), m(t), \alpha(t, x(t))) dt + \sigma d W(t).
\end{equation}
The corresponding mean field equilibrium is described as follows~\citep{lauriere2021numerical}.

\begin{definition}[MFG equilibrium]\label{def:MFG}
	The pair $ (\hat{\alpha}, \hat{m}) $ of the control input and the mean field is obtained as
	\begin{equation}\label{eqn:alpha_MFG}
	\hat{\alpha} = \arg\min_{\alpha} J^{MFG}_{\hat{m}}(\alpha),
	\end{equation}
	with the cost
	\begin{multline}\label{eqn:J_MFG}
	J^{MFG}_m(\alpha) = \E\bracket{\int_{0}^{T} f(x^{m}(t), m(t), \alpha(t, x^{m}(t))) dt + \right. \\
		\left. g(x^{m}(T), m(T)) },
	\end{multline} 
	where $ \hat{m}(t),\ t \in [0,T] $ is the distribution generated by $ x^{\hat{m}}(t) $ evolving according to the dynamics in \eqref{eqn:dx}.
\end{definition}

The solution is the mean field Nash equilibrium wherein:
\begin{itemize}
	\item Given the population distribution $ \hat{m} $, optimal control $ \hat{\alpha} $ for a representative agent is obtained in \eqref{eqn:alpha_MFG},
	\item The population distribution $ \hat{m} $ itself is generated by the optimal strategy $ \hat{\alpha} $.
\end{itemize}
Thus, this setup can be expressed as a fixed point problem. 

\subsection{Mean field control: social optimum}

Another problem can be formulated with a similar structure of the dynamics and the cost function. 
In the case of MFG, each agent optimizes their own cost resulting in a non-cooperative game. 
An alternative model arises when agents cooperate with each other to identify the best collective strategy, minimizing the average cost for the group. 

For MFG, the representative agent optimizes their cost in \eqref{eqn:J_MFG} assuming the strategies of other agents are fixed. 
In mean field control (MFC), however, the common feedback directly influences the population distribution through the mean field term.
Thus instead of a fixed point formulation, MFC can be cast as the solution of an optimal control problem driven by the stochastic dynamics. 

\begin{definition}[MFC equilibrium]\label{def:MFC}
	The optimal feedback control $ \alpha^* $ is determined as
	\begin{equation}\label{eqn:alpha_MFC}
	\alpha^* = \arg\min_{\alpha} J^{MFC}(\alpha),
	\end{equation}
	with the cost
	\begin{multline}\label{eqn:J_MFC}
	J^{MFC}(\alpha) = \E\bracket{\int_{0}^{T} f(x(t), m(t), \alpha(t, x(t))) dt + \right. \\
		\left. g(x(T), m(T)) },
	\end{multline} 
	where $ m(t),\ t \in [0,T] $ is the distribution generated by $ x(t) $ evolving according to the dynamics in \eqref{eqn:dx}.
	Furthermore, $ m^* $ is the mean-field corresponding to the optimal control $ \alpha^* $.
\end{definition}


\begin{definition}[Price of Anarchy]
	Let $ (\alpha^*, m^*) $ denote the pair of optimal control and the corresponding mean field obtained from MFC in \Cref{def:MFC}, and $ (\hat{\alpha}, \hat{m}) $ the MFG equilibrium from \Cref{def:MFG}. 
	It follows that
	\begin{equation}
	J^{MFC}(\alpha^*) = J^{MFG}_{m^*}(\alpha^*) \le J^{MFG}_{\hat{m}}(\hat{\alpha}),
	\end{equation}
	since the social optimum cost is lower for the representative agent.	
	The Price of Anarchy (PoA) is then defined as
	\begin{equation}\label{eqn:PoA}
	PoA = \frac{J^{MFG}_{\hat{m}}(\hat{\alpha})}{J^{MFC}(\alpha^*)} \ge 1,
	\end{equation}
	with similar structure of dynamics and costs.
\end{definition}

\subsection{PDE system for MFG}

Next, the optimality conditions for MFG in \Cref{def:MFG} can be expressed as a system of partial differential equations (PDE).

Define the Lagrangian as a function of the running cost:
\begin{equation}
L(x, m, \alpha, p) = f(x, m, \alpha) + \pair{b(x, m, \alpha), p},
\end{equation}
and the corresponding Hamiltonian for the representative agent is
\begin{equation}\label{eqn:H_transform}
H(x, m, p) = \max_{\alpha} \braces{-L(x, m, \alpha, p)}.
\end{equation}
It is assumed that $ H $ is strictly convex with respect to $ p $.

\begin{remark}\label{rem:PDE_MFG}
	In this setting of continuous state-action spaces, the optimality conditions are described by a set of coupled PDEs.
	\begin{align}
	-\frac{\partial u}{\partial t}(t, x) &- \nu \Delta u(t, x) + H(x, m(t, x), \nabla u(t, x)) = 0 \label{eqn:dudt_mfg} \\
	\frac{\partial m}{\partial t}(t, x) &- \nu \Delta m(t, x) \hfill \nonumber \\
	&- \diver \parenth{m(t, x) \partial_p H(x, m(t, x), \nabla u(t, x))} = 0 \label{eqn:dmdt_mfg} \\
	u(T, x) &= g(x, m(T, x)), \quad m(0, x) = m_0(x) \label{eqn:boun_mfg}
	\end{align}
	where $ \nu = \frac{\sigma^2}{2} $ and $ \partial_p H $ is the partial derivative of $ H $ w.r.t. the third variable which is $ p $ in \eqref{eqn:H_transform}. 
\end{remark}

Here, the optimal control for a typical agent satisfies the Hamilton-Jacobi-Bellman (HJB) equation through the value function, $ u(t,x) $.
Meanwhile, the density $ m(t,x) $ evolves according to the Fokker-Planck-Kolmogorov (FPK) equation. 
The HJB equation in \eqref{eqn:dudt_mfg} is formulated backward in time from a terminal value $ u(T, x) $, while the FPK equation in \eqref{eqn:dmdt_mfg} evolves forward in time from $ m(0, x) $.

Given the equilibrium values of population density and the value function, the optimal control for a representative agent is:
\begin{equation}\label{eqn:alpha_mfg}
\hat{\alpha}(t, x) = \arg\max_{\alpha} \braces{-L(x, m(t, x), \alpha, \nabla u(t, x))}.
\end{equation}

\subsection{PDE system for MFC}

Similarly, the optimality condition for MFC in \Cref{def:MFC} can also be expressed as a system of PDEs.
A key distinction in MFC is that when the control $ \alpha $ changes, the density $ m $ also varies within the optimization formulation of the representative agent in \eqref{eqn:J_MFC}.
Hence, the value function of the centralized optimizer explicitly depends upon the population distribution.

\begin{remark}\label{rem:PDE_MFC}
	The solution can be formulated as a modified set of PDEs:
	\begin{align}
	-\frac{\partial u}{\partial t}(t, x) &- \nu \Delta u(t, x) + H(x, m(t, x), \nabla u(t, x)) \nonumber \\
	& \int_\mathcal{D} \frac{\partial H}{\partial m} (\xi, m(t, x), \nabla u(t, \xi)) (x) m(t, \xi) d\xi = 0 \label{eqn:dudt_mfc}\\
	\frac{\partial m}{\partial t}(t, x) &- \nu \Delta m(t, x) \hfill \nonumber \\
	&- \diver \parenth{m(t, x) \partial_p H(x, m(t, x), \nabla u(t, x))} = 0 \label{eqn:dmdt_mfc} \\
	u(T, x) =& g(x, m(T, x)) + \int_\mathcal{D} \frac{\partial g}{\partial m} (\xi, m(T, x)) (x) m(T, \xi) d\xi \nonumber\\ 
	m(0, x) &= m_0(x) \label{eqn:boun_mfc}
	\end{align}
	where $ \nu = \frac{\sigma^2}{2} $.
\end{remark}

While the FPK equation \eqref{eqn:dmdt_mfc} retains the same structure, the HJB equation \eqref{eqn:dudt_mfc} has an additional term consisting of the derivative of $ H $ with respect to $ m $~\citep{bensoussan2013mean}. 
The optimal feedback control for a typical agent is then obtained as:
\begin{equation}\label{eqn:alpha_mfc}
\alpha^*(t, x) = \arg\max_{\alpha} \braces{-L(x, m(t, x), \alpha, \nabla u(t, x))}.
\end{equation}

As such, we consider problems with local mean-field interactions for which the solution to the MFG system \eqref{eqn:dudt_mfg}--\eqref{eqn:boun_mfg} exists and is unique (similarly for the MFC system \eqref{eqn:dudt_mfc}--\eqref{eqn:boun_mfc}).

\section{Mean Field Game and Control for Switching Hybrid Systems}\label{sec:ALG}

This section formulates mean field games for switching systems with pre-specified transitions, and it introduces a numerical framework to solve it. 
Finite-difference schemes are first constructed for the regular continuous flow, where the discrete value function and density are represented as grid functions on a uniform mesh. 
Discrete events introduce additional complexities into the coupled governing equations. 
To address these challenges, we design a specialized algorithm tailored to handle such scenarios effectively.  

\subsection{Switching Hybrid System}

Consider the system of governing equations for MFG, as given in \eqref{eqn:dudt_mfg},\eqref{eqn:dmdt_mfg}, which has a forward-backward structure with the initial and terminal conditions specified in \eqref{eqn:boun_mfg}. 
Now, suppose there are discrete events that may cause discontinuities in the dynamics.
These events may involve changes in spatial boundary conditions, alterations to the Hamiltonian's structure, or variations in the noise parameter $ \nu $, 
thereby causing sudden changes in the governing equation \eqref{eqn:dudt_mfg}, \eqref{eqn:dmdt_mfg} or the boundary conditions. 
We assume that these events occur at a prescribed set of time independent of the agents' individual states, and there is no reset of the state. 
The resulting system corresponds to a stochastic hybrid system without guard and reset, and it is referred to as switching hybrid system in this paper. 

Let the switching events occur at time indices $ \braces{e_1, \dots, e_{N_s}} $, where $ N_s $ denotes the total number of switches. 
The trajectory is divided into intervals between these switching events, with solutions computed separately for each interval and patched together across the entire time horizon. 
In short, the flow consists of a sequence of continuous flow sections concatenated with switching events.

The governing equations for each continuous flow section remain consistent, but are distinct from other sections. 
The solution is interconnected at the switching points by the initial and terminal conditions. 
Given the values of $ u(T, x) $ and $ m(0, x) $, additional constraints are imposed to match the trajectories of the value function and the density at the switching times. 
For example, $ u(t_{e_1}, x) $ and $ m(t_{e_1}, x) $ must align for the trajectory segments between $ (t_{e_0}, t_{e_1}) $ and $ (t_{e_1}, t_{e_2}) $, 
e.g., $u(t_{e_1}^-,x) = u(t_{e_1}^+, x)$. 

These matching conditions guarantee the consistency of the solution across switching events, allowing the solution to be concatenated seamlessly. 
In summary, we present a computational scheme to solve the forward-backward structure of the partial differential equations while addressing the switching dynamics and the additional internal boundary conditions. 

\subsection{Finite difference scheme for continuous flow}

We first present a computational scheme to solve the governing equation at each section of the continuous flow between switchings. 
For brevity, we do not introduce an additional notational convention to distinguish the variables over a specific section from other sections. 
It is understood that all of the developments in this subsection are devoted to a single section of continuous flow. 

We begin by considering 1-D spatial domain for derivation of the scheme.
This can be extended later to higher spatial dimensions.
Suppose the time duration of each section and the space variables are discretized into $ N_T + 1 $ and $ N_X + 1 $ points respectively.
Then the time step is $ \Delta t = T / N_T $ and the spacial grid size is $ h = (x_{max} - x_{min}) / N_X $.

The value function and the density are discretized by the matrices $ U, M \in \Re^{(N_T+1)\times(N_X+1)} $ as $ u(t_n, x_i) \approx U^n_i, m(t_n, x_i) \approx M^n_i $ for $ n \in \braces{0,\dots,N_T}, i \in \braces{0,\dots,N_X} $.

To evaluate the terms in the differential equations, we use the following relations:
\begin{itemize}
	\addtolength\itemsep{2mm}
	\item $ (D_t W)^n = (W^{n+1} - W^n)/ \Delta t, \quad W \in \Re^{N_T+1} $
	\item $ (D W)_i = (W_{i+1} - W_i)/ h, \quad W \in \Re^{N_X+1} $
	\item $ (\Delta_h W)_i = (W_{i+1} -2 W_i + W_{i-1})/ h^2, \quad W \in \Re^{N_X+1} $
	\item $ \bracket{\nabla_h W}_i = ((D W)_i, (D W)_{i-1}) $.
\end{itemize}

The discrete version of the Hamiltonian, which satisfies key numerical properties, is obtained from Godunov's scheme~\citep{laney1998computational}. 
For instance, the discrete Hamiltonian could have the following quadratic structure in 1-D:
\begin{equation*}
	\tilde{H}(x, m, \bracket{\nabla_h W}_i) = \tilde{H}(x, m,p_1, p_2) = \frac{1}{2} c_K ((p_1^-)^2 + (p_2^+)^2 ),
\end{equation*}
and its derivatives with respect to the $ p $-variables are 
\begin{equation*}
	\partial_{p_1} \tilde{H} = -c_K p_1^- ,\ \partial_{p_2} \tilde{H} = c_K p_2^+
\end{equation*}
where 
\begin{equation*}
	x^- = \max(-x, 0),\ x^+ = \max(x, 0).
\end{equation*}

Then the discrete version of the HJB equation~\eqref{eqn:dudt_mfg} for MFG can be formulated as~\citep{achdou2020mean}
\begin{align}\label{eqn:DU_mfg}
-(D_t U_i)^n &- \nu (\Delta_h U^n)_i + \tilde{H}(x_i, M^{n+1}, \bracket{\nabla_h U}_i) = 0,
\end{align}
$ \forall n \in \braces{0,\dots,N_{T}-1},\ i \in \braces{0,\dots,N_X} $ with the appropriate boundary conditions in space and the terminal condition $ U^{N_T}_i = g(x_i, M^{N_T}) $. 

Meanwhile, the discrete FPK equation~\eqref{eqn:dmdt_mfg} is expressed as
\begin{align}\label{eqn:DM_mfg}
(D_t M_i)^n &- \nu (\Delta_h M^{n+1})_i - \mathcal{T}_i(U^n, M^{n+1}) = 0,
\end{align}
$ \forall n \in \braces{0,\dots,N_{T}-1},\ i \in \braces{0,\dots,N_X} $ starting from $ M^0_i = m_0(x_i) $.
Here, the discrete transport operator $ \mathcal{T}_i $ has the following structure
\begin{multline}\label{eqn:Ti}
\mathcal{T}_i(U, M) = \frac{1}{h}\parenth{M_i \partial_{p_1}\tilde{H}(x_i, M_i, \bracket{\nabla_h U}_i ) \right. \\ 
	\left.  -M_{i-1} \partial_{p_1}\tilde{H}(x_{i-1}, M_{i-1}, \bracket{\nabla_h U}_{i-1} ) \right.\\
	\left.  +M_{i+1} \partial_{p_2}\tilde{H}(x_{i+1}, M_{i+1}, \bracket{\nabla_h U}_{i+1} ) \right.\\		
	\left.  -M_{i} \partial_{p_2}\tilde{H}(x_{i}, M_{i}, \bracket{\nabla_h U}_{i} ) 
},
\end{multline}
which ensures that the total mass of the density matrix is constant over time.

Both the discrete HJB and FPK equations are implicit schemes in time.
When \eqref{eqn:DU_mfg}, \eqref{eqn:DM_mfg} are combined, the system of equations can be expressed as a root finding problem:
\begin{equation}\label{eqn:root_Z}
\varphi(Z) = 0,
\end{equation}
where $ Z = (U, M) $ and $ \varphi : \Re^{2\times(N_T+1)\times(N_X+1)} \to \Re^{2\times(N_T+1)\times(N_X+1)} $.

The simplest method to obtain the solution is the Picard iteration approach, where each of the equations is iteratively updated using the most recent estimates of the value function and the density. 
However, this procedure has significant challenges in terms of convergence to a solution. 
Therefore, the chosen robust approach is to solve via Newton iterations from an initial guess.
Starting from $ Z^{(0)} $ where $ Z^{(k)} = (U^{(k)}, M^{(k)}) $, the solution is the asymptotic iterate of
\begin{equation}\label{eqn:newton}
Z^{(k+1)} = Z^{(k)} - J_{\varphi}(Z^{(k)})^{-1} \varphi(Z^{(k)}).
\end{equation}
where $ J_{\varphi}(Z^{(k)}) $ is the Jacobian of $ \varphi $ at $ Z^{(k)} $.

\subsection{Computational algorithm for switching}

\begin{algorithm}[h!]
	\linespread{1.25}\selectfont
	\Input Initial and terminal condition : $ U^{N_T}, M^0 $
	\begin{algorithmic}[1]
		\State Initial guess $ U^{(0)}, M^{(0)} $
		\For {$ k = 0$ \textbf{to} $N_{iter}-1$}
		\FnBlock{$ \varPhi (Z) $}
		\State $ U^{e_{N_s+1}} \gets U^{N_T},\ M^{e_0} \gets M^0 $
		\For {$ s = 0$ \textbf{to} $N_s$}
		\State $ t_s \gets \braces{t_{e_s}, \dots, t_{e_{s+1}}} $
		\State $ U^{e_{s+1}} = U[t_{e_{s+1}}],\ M^{e_s} = M[t_{e_s}] $
		\State $ X[t_s] \gets \varphi_s\parenth{(U[t_s], M[t_s])} $
		\EndFor
		\EndFnBlock{ X}
		\State $ \Delta Z^{(k)} \gets \text{solve} \braces{ J_{\varPhi}(Z^{(k)}) \ Y = \varPhi(Z^{(k)}) } $ for $ Y $
		\State $ Z^{(k+1)} \gets Z^{(k)} - \Delta Z^{(k)} $
		\EndFor
		\State \Return $ (U^{(N_{iter})}, M^{(N_{iter})}) = Z^{(N_{iter})} $
	\end{algorithmic}
	\caption{Solver for switching systems}
	\label{alg:sim}
\end{algorithm}

Now we integrate the above numerical solution with an algorithm to address the switching events. 
The proposed framework is summarized in \Cref{alg:sim}, where the solution to the full system is expressed as the root of $ \varPhi(Z) = 0 $.
The expanded function $ \varPhi $ concurrently computes the value of each individual term in \eqref{eqn:root_Z} for every partition $ (t_{e_s}, t_{e_{s+1}}) $, under the constraint that initial and terminal conditions of $ U^{e_{s+1}}, M^{e_s} $ match. 
The domain of the one-to-one function $ \varPhi $ is $ \Re^{2\times(N_T+1)\times(N_X+1)} $, whereas the domain of $ \varphi_s $ is $ \Re^{2\times(N_{t_s}+1)\times(N_X+1)} $, where $ N_{t_s} $ represents the number of time steps in each continuous flow section. 

In numerical implementation, evaluating the inverse of the Jacobian of $ \varPhi $ is computationally expensive and unnecessary.
Instead, the solution in Step 11 of \Cref{alg:sim} can be computed by solving the linear system of equations
\begin{equation}\label{eqn:jac_vec}
J_{\varPhi}(Z^{(k)}) Y = \varPhi(Z^{(k)}),
\end{equation}
to obtain $ Y = J_{\varPhi}(Z^{(k)})^{-1} \varPhi(Z^{(k)}) $. 
Moreover, iterative solvers are available for this linear system that only require the Jacobian-vector product, thereby avoiding the need to construct the full Jacobian matrix explicitly. 

The choice of the initial value $ Z^{(0)} $ is crucial in the first step.
Since the solution converges efficiently for large values of viscosity $ \nu $, the system is first solved with a large $ \nu $, which then serves as the initial guess for the case with smaller $ \nu $.
Thus a continuation method is applied to gradually reach to the required level of $ \nu $. 

\section{Numerical Example}\label{sec:SIM}

In this final section, we formulate the problem of emergency evacuation from a room, accounting for congestion effects and transitions in the boundary conditions.
Consider the example of a crowd exiting a room with obstacles, such as a movie theater during a fire emergency.
The problem is initially formulated for the smooth case without any transition events. 

\begin{figure}[h!]
	\centering
	\subfigure[Initial geometry]{
		\includegraphics[height=0.425\columnwidth, width=0.45\linewidth]{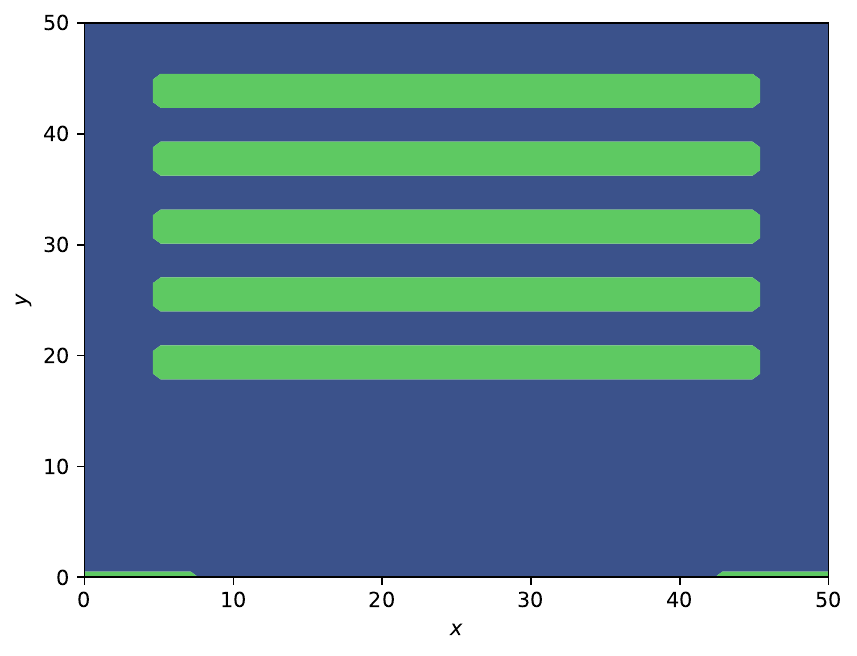}}
	\hfill
	\subfigure[Initial density]{
		\includegraphics[height=0.425\columnwidth, width=0.5\linewidth]{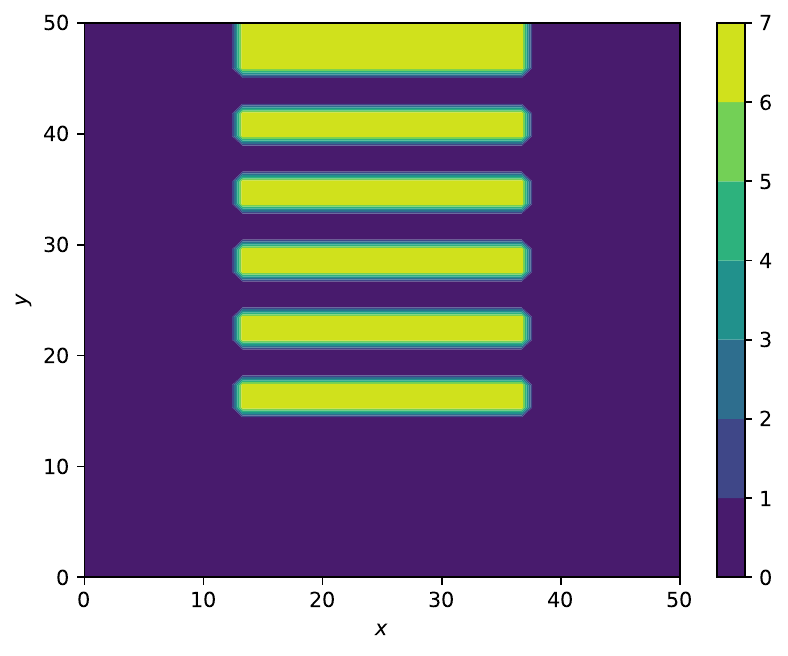}}
	\caption{Room evacuation with obstacles and two exits}
	\label{fig:init_conf}
\end{figure}

Let the geometry of the room be a square of the size of $ 50 $m on each side, featuring two doors which serve as exits. 
These doors are initially positioned along the bottom side at $(x,y) \in [0, 7.5] \cup [42.5, 50]\times \{0\} $.
Furthermore, the room also contains obstacles (represented in green in \Cref{fig:init_conf}.(a)) which are impassable areas.
Meanwhile, the initial distribution of people at $ t = 0 $ is uniform in the space adjoining these obstacles as shown in \Cref{fig:init_conf}.(b).
That is, they are situated with equal density in the yellow areas, while the density is zero otherwise. 

For this illustration, we assume that the drift term in the dynamics is equal to the control, i.e., $ b = \alpha $ in \eqref{eqn:dx}.
The running cost in the optimization problem of \eqref{eqn:J_MFG} has the structure:
\begin{equation}\label{eqn:f_evac}
f(x, m, \alpha) = \frac{\norm{\alpha}^2}{32}(1+m)^{3/4} + \frac{1}{3200}.
\end{equation}
So the cost increases with both the control magnitude and the local density. 
Thus congestion effects arise from local interaction with the mean field.
There is also a fixed term in \eqref{eqn:f_evac} that penalizes the trivial strategy of remaining stationary.
Hence, each individual wants to exit the room in the minimum possible time. 
Finally, there is no terminal cost in the system, $ g \equiv 0 $.

This setup enables the study of how people will evacuate optimally, given the constraints, and how congestion and obstacles affect the evacuation dynamics. 
In the limit of a large number of people, this problem can be formulated as a system of coupled partial differential equations (PDEs). 
In this case, the Hamiltonian is evaluated from \eqref{eqn:H_transform} as
\begin{equation}
H(x, m, p) = \frac{8 \norm{p}^2}{(1+m)^{3/4}} - \frac{1}{3200}.
\end{equation}

Thus, the governing equations for the MFG system from \eqref{eqn:dudt_mfg},\eqref{eqn:dmdt_mfg} can be written as
\begin{align}
-\frac{\partial u}{\partial t} &- \nu \Delta u + \frac{8 \norm{\nabla u}^2}{(1+m)^{3/4}} = \frac{1}{3200} \label{eqn:dudt_evac} \\
\frac{\partial m}{\partial t} &- \nu \Delta m -\diver \parenth{\frac{16 m \nabla u}{(1+m)^{3/4}}} = 0, \label{eqn:dmdt_evac}
\end{align}
while for MFC, the HJB equation from \eqref{eqn:dudt_mfc} is :
\begin{align}
-\frac{\partial u}{\partial t} - \nu \Delta u + \parenth{\frac{8}{(1+m)^{3/4}} - \frac{6m}{(1+m)^{7/4}}} \norm{\nabla u}^2 = \frac{1}{3200} \label{eqn:dudt_emfc}
\end{align}
with a value of $ \nu = 0.05 $.

At the exit doors of the room, we apply Dirichlet boundary condition of $ u = 0 $ for the value function, and $ m = 0 $ since the density outside is assumed to be zero.
Meanwhile, at the walls of the environment and the rectangular obstacles, a Neumann condition is imposed on both $ u, m $, ensuring that the velocity of the agents remains tangential to the surface.
As both $ \deriv{u}{n}, \deriv{m}{n} = 0 $, there is no flux of density through the edges.

The simulation time is set to $ T = 50 $ minutes, with the terminal condition for the value function being $ 0 $ since there is no terminal cost.
However, the problem is now extended to incorporate certain events such as changes in available exits due to the unfolding emergency.

\begin{figure}[h!]
	\centering
	\includegraphics[width=\linewidth]{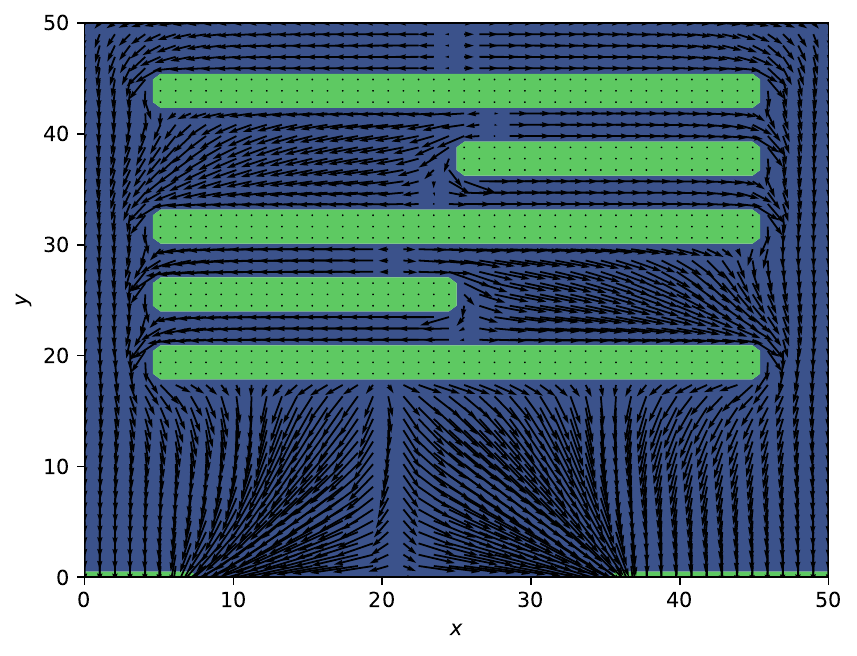}
	\caption{Geometry and velocity plot at $ t = 10 $ minutes}
	\label{fig:vel_plot}
\end{figure}

Specifically, we consider two events which occur at $ \braces{t_{e1}, t_{e2}} $.
At $ t_{e1} $, two parts of the obstacles shown in \Cref{fig:init_conf}.(a) are suddenly removed from the environment.
Later, at $ t_{e2} $, the exit door on the right hand side of the system is widened to double its original size (see \Cref{fig:vel_plot} for the updated geometry after both events).
These events change the Neumann and Dirichlet boundary conditions of the differential equations respectively.
So they will need to addressed systematically to obtain the full solution.
The switching times are $ t_{e1} = 2,\ t_{e2} = 5 $ in minutes.

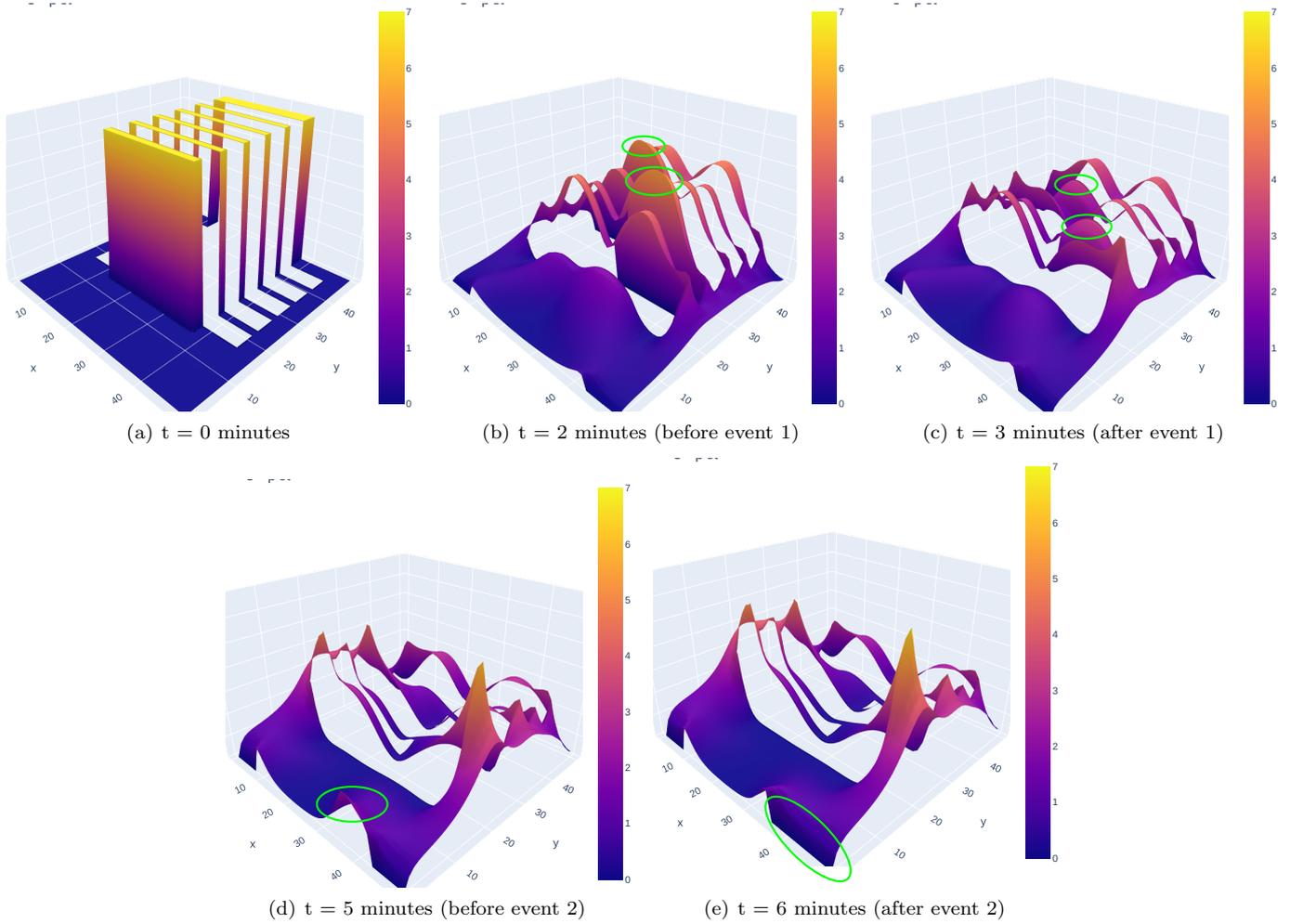
\begin{figure*}[h!]
	\centering
	\subfigure[t = 0 minutes]{
		\includegraphics[width=0.32\linewidth]{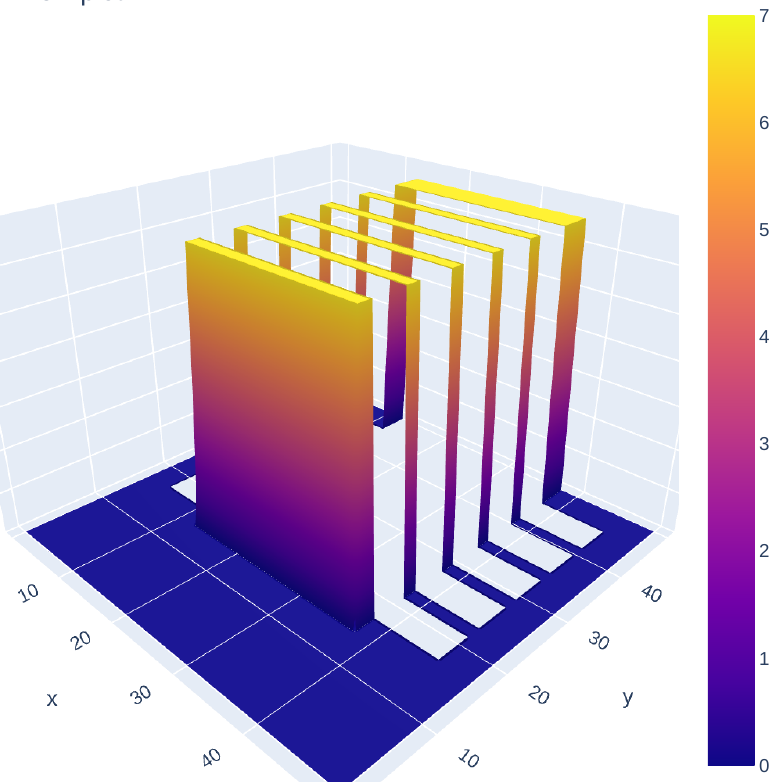}}
	\hfill
	\subfigure[t = 2 minutes (before event 1)]{
        \begin{tikzpicture}
            \node[anchor=south west, inner sep=0] (figb) at (0,0) {\includegraphics[width=0.32\linewidth]{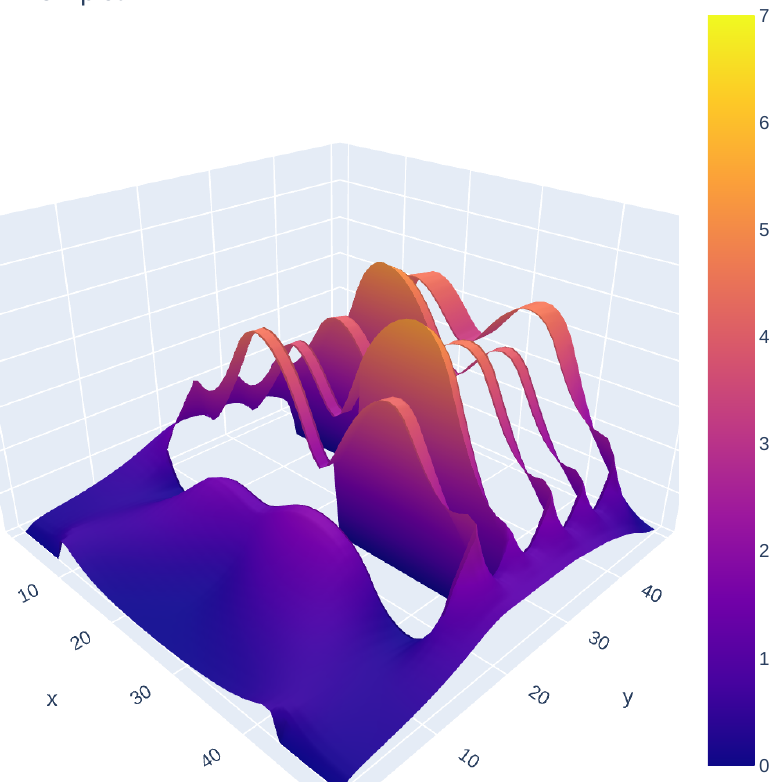}};
            \draw[green, thick] (3.05,3.3) ellipse (0.4cm and 0.2cm);
            \draw[green, thick] (2.9,3.8) ellipse (0.3cm and 0.14cm);
    \end{tikzpicture}}
	\hfill
	\subfigure[t = 3 minutes (after event 1)]{
		\begin{tikzpicture}
		\node[anchor=south west, inner sep=0] (figc) at (0,0) {\includegraphics[width=0.32\linewidth]{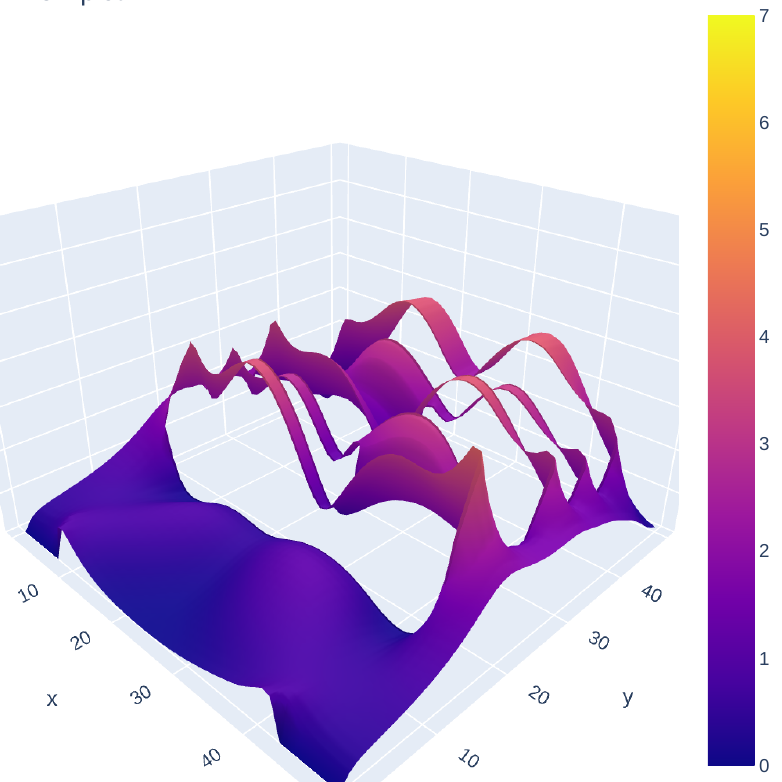}};
		\draw[green, thick] (3.05,2.65) ellipse (0.35cm and 0.17cm);
		\draw[green, thick] (2.9,3.25) ellipse (0.3cm and 0.14cm);
		\end{tikzpicture}}
	%
	\\
	\centering
	\subfigure[t = 5 minutes (before event 2)]{
		\begin{tikzpicture}
		\node[anchor=south west, inner sep=0] (figd) at (0,0) {\includegraphics[width=0.32\linewidth]{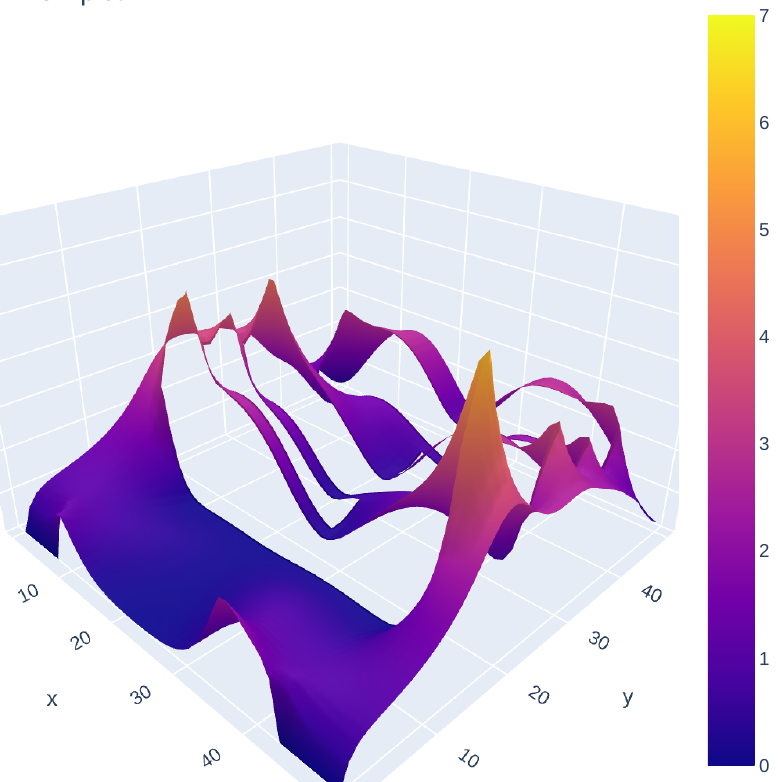}};
		\draw[green, thick] (1.8,1.2) ellipse (0.5cm and 0.25cm);
		\end{tikzpicture}}
	\subfigure[t = 6 minutes (after event 2)]{
		\begin{tikzpicture}
		\node[anchor=south west, inner sep=0] (figd) at (0,0) {\includegraphics[width=0.32\linewidth]{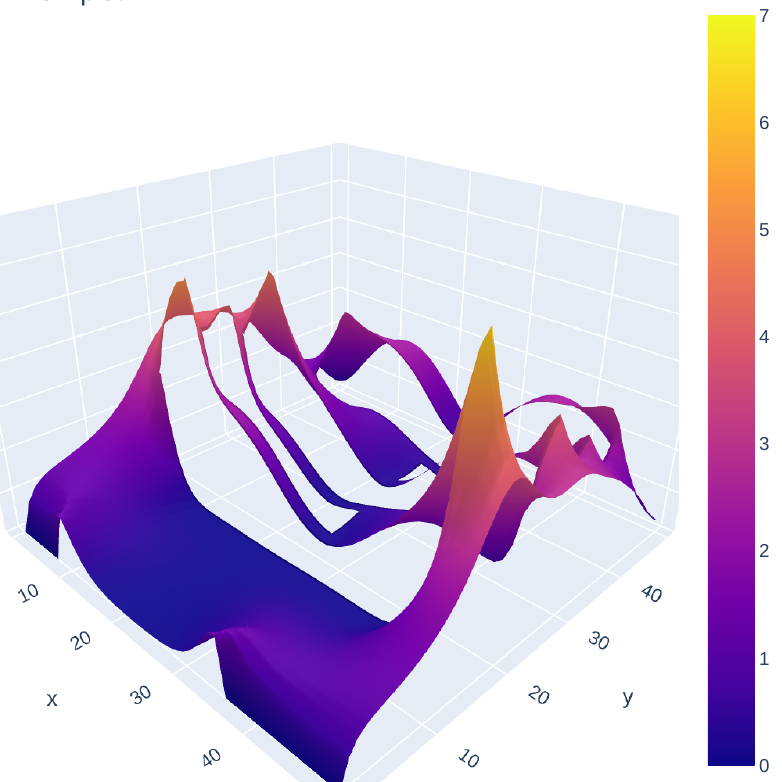}};
		\draw[green, thick, rotate around={-45:(2.2,0.4)}] (2.2,0.4) ellipse (0.8cm and 0.3cm);
		\end{tikzpicture}}
	\caption{Density evolution across the events}
	\label{fig:density}
\end{figure*}

\Cref{alg:sim} is implemented to obtain the solution for the dynamics in \eqref{eqn:dudt_evac},\eqref{eqn:dmdt_evac} with these transitions.
However, the Newton iterations require a good initial guess $ Z^{(0)} $ for the scheme to work properly. 
Hence, the system is split into separate problems based on the corresponding intervals $ \braces{ [0, t_{e_1}], [t_{e_1}, t_{e_2}], [t_{e_2}, T] }$, and their individual solutions are patched together to form the initial guess for the full system.

The velocity plot of the solution is shown in \Cref{fig:vel_plot} at a time after both events have occurred. 
It is observed that people move between the obstacles towards the corridors at the sides and then use the doors.
The density evolution at certain points in time is shown in \Cref{fig:density}.
It is found that the scheme effectively captures congestion effects and ensures that the velocity remains tangential to the walls.

The switching nature of the system introduces very interesting properties in the final solution. 
People are rational agents and they have choice of deciding which path to take in order to minimize the cost.
So it is expected that they will take advantage of the fact that some obstacles are removed in the first event at $ t_{e_1} $.
This is exactly what we observe in \Cref{fig:density}.(b),(c), where people gather above the area where the obstacle is expected to be removed, as represented by the area of higher density marked by green ovals.
Thus, they can move out quickly after the event happens. 

Similarly, a higher density is seen near the right-side door in \Cref{fig:density}.(d) before the second event, since that door is going to be widened to $(x,y) \in [35, 50]\times \{0\} $. 
Furthermore, the congestion effects manifest in the form of lower velocity in areas where density is higher. 
This is represented by the length of arrows in \Cref{fig:vel_plot}, for example, the length is smaller in the intersection of the corridor with the path between the obstacles. 

\begin{figure}[h!]
	\centering
	\includegraphics[width=\linewidth]{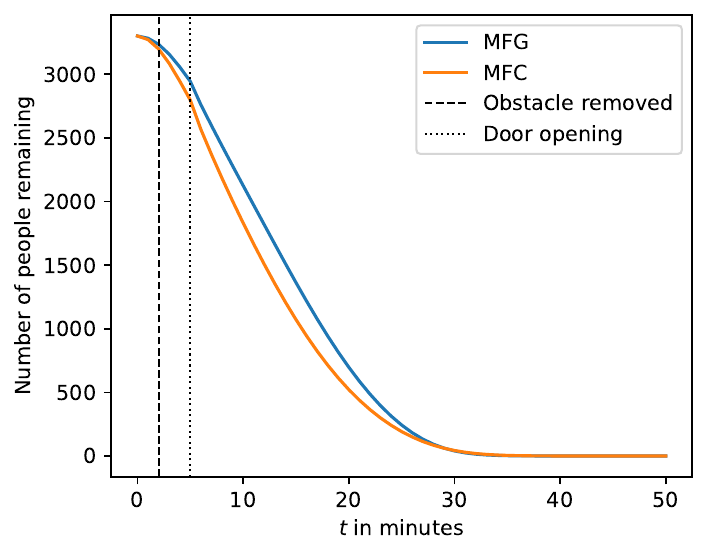}
	\caption{Comparison of total number of people remaining}
	\label{fig:compare_density}
\end{figure}

Meanwhile, \Cref{alg:sim} can be used to obtain the solution for the mean field control problem in \eqref{eqn:dudt_emfc} with the same boundary condition transitions. 
This scenario corresponds to a social optimum, where individuals coordinate with each other to determine the best possible strategy. 
Initially, the total number of people in the room is $ 3300 $.
It is numerically verified that MFC leads to a faster and more efficient evacuation of the crowd in comparison to MFG (see \Cref{fig:compare_density}).
In this case, the Price of Anarchy calculated  using \eqref{eqn:PoA} turns out to be PoA $ = 1.0181 $.
Additionally, when the second event occurs and the right door is widened, a discontinuity appears in the slope of the curve, and the evacuation rate increases. 

\begin{figure}[h!]
	\centering
	\includegraphics[width=\linewidth]{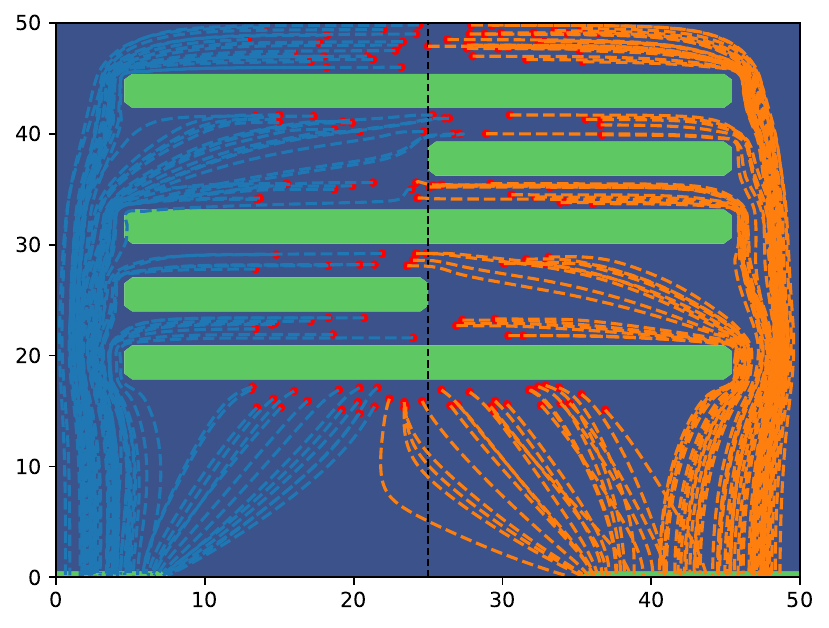}
	\caption{Sample trajectories (color coded by the exit taken)}
	\label{fig:sample_traj}
\end{figure}

Once the full solution for the density and the value function has been obtained, individual trajectories can be analyzed to understand their evolution.
The control which is given by $ b = - H_p(x, m, \nabla u) $ is numerically evaluated at each point along the trajectory using the solutions for $ m, u $.
Multiple initial conditions are considered with the assumption of zero noise for numerical purposes, and their respective paths are shown in \Cref{fig:sample_traj}.
As observed, a larger proportion of people prefer the right side exit since that door is widened during the second event.
Furthermore, individuals also tend to choose slightly longer paths to avoid congestion, particularly when they are positioned near the area where obstacles are removed during the first event. 

\section{Conclusion}

This paper discussed the mean field game framework for switching hybrid systems with pre-specified transitions.
We proposed an algorithm to obtain scalable control strategies that enable agents to work together, adjusting their behavior based on overall objectives and environmental conditions.
The capability of the computational methodology was demonstrated through a numerical illustration of room evacuation, where changes in the boundary conditions led to distinctive characteristics in crowd evolution. 
Future work includes extension to general hybrid systems, for example, when each agent can switch between discrete modes (e.g., ``stressed" vs ``calm" strategies in this example). 

\bibliography{references}

\end{document}